\documentclass{article}
\usepackage{array,amsmath}
\numberwithin{equation}{section}
\usepackage{indentfirst}

\begin{document}

\newtheorem{thm}{Theorem}[section]
\newtheorem{cor}[thm]{Corollary}
\newtheorem{prop}[thm]{Proposition}
\newtheorem{conj}[thm]{Conjecture}
\newtheorem{lem}[thm]{Lemma}
\newtheorem{Def}[thm]{Definition}
\newtheorem{rem}[thm]{Remark}
\newtheorem{prob}[thm]{Problem}
\newtheorem{ex}{Example}[section]

\newcommand{\be}{\begin{equation}}
\newcommand{\ee}{\end{equation}}
\newcommand{\ben}{\begin{enumerate}}
\newcommand{\een}{\end{enumerate}}
\newcommand{\beq}{\begin{eqnarray}}
\newcommand{\eeq}{\end{eqnarray}}
\newcommand{\beqn}{\begin{eqnarray*}}
\newcommand{\eeqn}{\end{eqnarray*}}
\newcommand{\bei}{\begin{itemize}}
\newcommand{\eei}{\end{itemize}}

\newcommand{\pa}{{\partial}}
\newcommand{\V}{{\rm V}}
\newcommand{\R}{{\bf R}}
\newcommand{\K}{{\rm K}}
\newcommand{\e}{{\epsilon}}
\newcommand{\tomega}{\tilde{\omega}}
\newcommand{\tOmega}{\tilde{Omega}}
\newcommand{\tR}{\tilde{R}}
\newcommand{\tB}{\tilde{B}}
\newcommand{\tGamma}{\tilde{\Gamma}}
\newcommand{\fa}{f_{\alpha}}
\newcommand{\fb}{f_{\beta}}
\newcommand{\faa}{f_{\alpha\alpha}}
\newcommand{\faaa}{f_{\alpha\alpha\alpha}}
\newcommand{\fab}{f_{\alpha\beta}}
\newcommand{\fabb}{f_{\alpha\beta\beta}}
\newcommand{\fbb}{f_{\beta\beta}}
\newcommand{\fbbb}{f_{\beta\beta\beta}}
\newcommand{\faab}{f_{\alpha\alpha\beta}}

\newcommand{\pxi}{ {\pa \over \pa x^i}}
\newcommand{\pxj}{ {\pa \over \pa x^j}}
\newcommand{\pxk}{ {\pa \over \pa x^k}}
\newcommand{\pyi}{ {\pa \over \pa y^i}}
\newcommand{\pyj}{ {\pa \over \pa y^j}}
\newcommand{\pyk}{ {\pa \over \pa y^k}}
\newcommand{\dxi}{{\delta \over \delta x^i}}
\newcommand{\dxj}{{\delta \over \delta x^j}}
\newcommand{\dxk}{{\delta \over \delta x^k}}

\newcommand{\px}{{\pa \over \pa x}}
\newcommand{\py}{{\pa \over \pa y}}
\newcommand{\pt}{{\pa \over \pa t}}
\newcommand{\ps}{{\pa \over \pa s}}
\newcommand{\pvi}{{\pa \over \pa v^i}}
\newcommand{\ty}{\tilde{y}}
\newcommand{\bGamma}{\bar{\Gamma}}

\font\BBb=msbm10 at 12pt
\newcommand{\Bbb}[1]{\mbox{\BBb #1}}

\newcommand{\qed}{\hspace*{\fill}Q.E.D.}  

\title{Some important applications of improved Bochner inequality on Finsler manifolds}
\author{ Xinyue Cheng \footnote{supported by the National Natural Science Foundation of China (11871126) and the Science Foundation of Chongqing Normal University (17XLB022)}}
\date{}

\maketitle

\begin{abstract}
We establish some important inequalities under the condition that the weighted Ricci curvature $\mathrm{Ric}_{\infty}\geq K$ for some constant $K >0$ by using improved Bochner inequality and its integrated form. Firstly, we obtain a sharp Poincar\'{e}-Lichnerowicz inequality. Further, we give a new proof for logarithmic Sobolev inequality. Finally, we obtain an estimate of the volume of geodesic balls. 
{\bf MR(2010) Subject Classification:}  53B40, 53C60
\end{abstract}

\section{Introduction}

In Riemann geometry, Bochner formula is a bridge to use the analytic tools to study the geometry and topology of a manifold, which is stated as follows: for any smooth function $f$, we have
\be
\frac{1}{2}\Delta |\nabla f|^{2}= |{\rm Hess}  f|^{2}+\langle \nabla \Delta f, \nabla f\rangle +Ric(\nabla f, \nabla f),
\ee
where Hess$f$ is the Hessian of $f$ and $\Delta f$ is the Laplacian of $f$. Making use of the Bochner formula to harmonic 1-form, one can show that, if a compact $n$-dimensional Riemannian manifold  has ${\rm Ric}\geq 0$, then $b_{1}(M)\leq n$ and ``= " holds if and only if $M$ is isometric to a flat torus. Bochner formula also gives an analytic approach to define the Ricci curvature.

Further, the Bochner-Weitzenb\"{o}ck formula and the corresponding Bochner inequality on Finsler manifolds were established by Ohta-Sturm in \cite{OS2}. Let $(M, F, m)$ be an $n$-dimensional Finsler manifold equipped with a smooth measure $m$. Given $u\in {\cal C}^{\infty}(M)$, Ohta-Sturm obtained the following point-wise Bochner-Weitzenb\"{o}ck formula
\be
\Delta^{\nabla u}\left[\frac{F^{2}(\boldsymbol{\nabla} u)}{2}\right]-d(\boldsymbol{\Delta} u)(\boldsymbol{\nabla} u)=\operatorname{Ric}_{\infty}(\boldsymbol{\nabla} u)+\left\|\boldsymbol{\nabla}^{2} u\right\|_{\mathrm{HS}(\boldsymbol{\nabla} u)}^{2}
\ee
as well as
\be
\Delta^{\nabla u}\left[\frac{F^{2}(\boldsymbol{\nabla} u)}{2}\right]-d(\Delta u)(\boldsymbol{\nabla} u) \geq \operatorname{Ric}_{N}(\boldsymbol{\nabla} u)+\frac{(\Delta u)^{2}}{N}
\ee
for $N \in(-\infty, 0) \cup[n, \infty],$  which are both pointwise on $M_{u}:=\{x\in M \ | du (x)\neq 0\}$. Here $\|\cdot\|_{H S(\nabla u)}$ stands for the Hilbert-Schmidt norm with respect to $g_{\nabla u}$ and they used the nonlinear Laplacian $\Delta$ and its linearization $\Delta^{\nabla u}.$ Then Ohta-Sturm shown the integrated form of the Bochner-Weitzenb\"{o}ck formula and the corresponding Bochner inequality (\cite{OS2}). Furthermore, one can get the following improved Bochner inequality.

\begin{lem}\label{Pro12.13} {\rm (Improved Bochner inequality \cite{Oh3})}   Assume that ${\rm Ric}_{\infty} \geq K$ for some $K \in {R}.$  Then, for any $u \in {C}^{\infty}(M)$, we have
\be
\Delta^{\nabla u}\left[\frac{F^{2}({\nabla} u)}{2}\right]-d({\Delta} u)({\nabla} u) \geq K F^{2}({\nabla} u)+d[F({\nabla} u)]\left(\nabla^{\nabla u}[F({\nabla} u)]\right)  \label{eq12.13}
\ee
pointwise on $M_{u}$.
\end{lem}

Based on this, one can obtain the integrated form of the improved inequality as follows.

\begin{lem}\label{IntformBo} {\rm (Integrated form \cite{Oh3})} Assume ${\rm Ric}_{\infty} \geq K$ for some $K \in {R}.$ Then, given $u \in H_{\mathrm{loc}}^{2}(M) \cap {\cal C}^{1}(M)$ such that $\Delta u \in H_{\mathrm{loc}}^{1}(M),$ we have
\beq
&&-\int_{M} d \phi\left(\nabla^{\nabla u}\left[\frac{F^{2}({\nabla} u)}{2}\right]\right) d{m} \nonumber\\
&& \geq \int_{M} \phi\left\{d({\Delta} u)({\nabla} u)+K F^{2}({\nabla}u)+d[F({\nabla}u)]\left(\nabla^{\nabla u}[F({\nabla} u)]\right)\right\} d{m} \label{Intform}
\eeq
for all nonnegative functions $\phi \in H_{c}^{1}(M) \cap L^{\infty}(M)$. Here, $\phi$ is called the test function.
\end{lem}

The Bochner-Weitzenb\"{o}ck formula and the corresponding Bochner inequality on Finsler manifolds have been applied to many important research topics. For example, following Bochner-Weitzenb\"{o}ck type formula, Wang-Xia give a sharp lower bound for the first (nonzero) Neumann eigenvalue of Finsler-Laplacian in Finsler manifolds in terms of diameter, dimension, weighted Ricci curvature (\cite{WangXia}). Q. Xia gives an (integrated) $p(>1)$-Bochner-Weitzenb\"{o}ck formula and the $p$-Reilly type formula on Finsler manifolds. As applications, Q. Xia obtains the $p$ Poincar\'{e} inequality on an $n$-dimensional compact Finsler manifold without boundary or with convex boundary under the assumption that ${\rm Ric}_{N} \geq K$ for $N \in[n, \infty]$ and $K \in {R}$ (see \cite{Xia}).

Lemma \ref{Pro12.13} and Lemma \ref{IntformBo} are  important foundations for our discussions in this paper. Starting from (\ref{Intform}), we can get the following inequality firstly.

\begin{thm}\label{PLineq} {\rm (Poincar\'{e}-Lichnerowicz inequality)} Suppose that $(M, F, m)$ is compact and satisfies $m(M)=1$ and $\operatorname{Ric}_{\infty} \geq K>0$. Then, for any $f \in H^{1}(M)$, we have the following
\be
\operatorname{Var}_{m}(f) \leq \frac{1}{K} \int_{M} F^{2}({\nabla} f) d{m}-\frac{2}{K}\int_{M}\left(\int_{0}^{\infty} g(t)dt\right)dm, \label{PoLiin}
\ee
where ${\rm Var}_{m}(f)$ denotes the variance of $f$ and
\[
g(t):= g_{\nabla u_{t}}\left(\nabla ^{{\nabla u_{t}}}F({\nabla u_{t}}), \nabla ^{{\nabla u_{t}}}F({\nabla u_{t}})\right),
\]
$\left(u_{t}\right)_{t \geq 0}$ is the global solution to the heat equation with $u_{0}=f$.
\end{thm}

In 2017, Ohta proved the following proposition (see Proposition 4.1 in \cite{Oh2}): Assume that $M$ is compact, $\operatorname{Ric}_{\infty} \geq K>0,$ and
\be
\int_{M} \frac{F^{2}({\nabla} u)}{u} d {m} \leq-C \int_{M}\left\{d u\left(\nabla^{\nabla u}\left[\frac{F^{2}({\nabla} u)}{2 u^{2}}\right]\right)+u \cdot d[{\Delta}(\log u)]({\nabla}[\log u])\right\} d{m} \label{LSi1}
\ee
holds for some constant $C>0$ and all functions $u \in H^{2}(M) \cap \mathcal{C}^{1}(M)$ such that $\Delta u \in H^{1}(M)$ and $\inf _{M} u>0 .$ Then the logarithmic Sobolev inequality
\be
\int_{\{f>0\}} f \log f d{m} \leq \frac{C}{2} \int_{\{f>0\}} \frac{F^{2}({\nabla} f)}{f} d{m} \label{LSi2}
\ee
holds for all nonnegative functions $f \in H^{1}(M)$ with $\int_{M} f d{m}=1$. Further, based on the above result, we can give a new proof for the following  logarithmic Sobolev inequality under the condition that ${\rm Ric}_{\infty} \geq K>0$.

\begin{thm}\label{LSIne}{\rm (Logarithmic Sobolev inequality)} Assume that $(M, F, {m})$ is compact and satisfies $\mathrm{Ric}_{\infty} \geq K>0$ and ${m}(M)=1$. Then we have
\be
\int_{M} f \log f d{m} \leq \frac{1}{2 K } \int_{M} \frac{F^{2}({\nabla} f)}{f} d{m} \label{eq16.7m}
\ee
for all nonnegative functions $f \in H^{1}(M)$ with $\int_{M} f d{m}=1$.
\end{thm}

Fix $p\in M$, let $r(x):= d_{F}(p, x)$ be the distance function on $M$ from $p$. It is well-known that $r$ is smooth on $M\backslash \{p\}$ away from the cut points of $p$. Let $B_{p}(R)$ be the forward geodesic ball of $M$ with radius $R$ centered at $p .$ The volume of $B_{p}(R)$ with respect to $m$ is defined by
\[
{\rm vol}\left(B_{p}(R)\right)=\int_{B_{p}(R)} d m .
\]
We have the following estimate of the volume of geodesic ball $B_{p}(R)$.
\begin{thm}\label{volgeoBall}   Let $(M, F, m)$ be a positively complete Finsler manifold of dimension $n$ and $r (x)=d_{F}(p, x)$ be the distance
function from a fixed point p. Suppose that  ${\rm Ric}_{\infty} \geq K>0$.  Then
\be
{\rm vol}\left(B_{p}(R)\right)\leq \frac{1}{K}\| \Delta r\|^{2}_{L^{2}(B_{p}(R))}. \label{geoBall}
\ee
\end{thm}

\vskip 3mm

The paper is organized as follows. In Section \ref{Sec2}, we give some necessary definitions and notations. Then the proof of a sharp Poincar\'{e}-Lichnerowicz inequality is given in Section \ref{Sec3}. Further, we give a new proof for logarithmic Sobolev inequality in Section \ref{Sec4} which is different from 
Ohta's proof by the curvature-dimension condition in \cite{Oh1}. Finally, an estimate of the volume of geodesic balls is obtained in Section \ref{Sec5}.

\section{Preliminaries}\label{Sec2}

Let $M$ be an $n$-dimensional manifold. A Finsler metric $F$ on $M$ is a non-negative function on $TM$ such that $F$ is ${\cal C}^{\infty}$ on $TM\backslash \{0\}$ and the restriction $F_{x}:=F|_{T_{x}M}$ is a Minkowski function on $T_{x}M$ for all $x\in M$. For Finsler metric $F$ on $M$, there is a Finsler co-metric $F^{*}$ on $M$ which is non-negative function on the cotangent bundle $T^{*}M$ given by
\be
F^{*}(x, \xi):=\sup\limits_{y\in T_{x}M\setminus \{0\}} \frac{\xi (y)}{F(x,y)}, \ \ \forall \xi \in T^{*}_{x}M.
\ee
We call $F^{*}$ the dual Finsler metric of $F$. Finsler metric $F$ and its dual Finsler metric $F^{*}$ satisfy the following relation.

\begin{lem}{\rm (Lemma 3.1.1, \cite{shen1})}\label{shen311} Let $F$ be a Finsler metric on $M$ and $F^{*}$ its dual Finsler metric. For any vector $y\in T_{x}M\setminus \{0\}$, $x\in M$, the covector $\xi =g_{y}(y, \cdot)\in T^{*}_{x}M$ satisfies
\be
F(x,y)=F^{*}(x, \xi)=\frac{\xi (y)}{F(x,y)}. \label{shenF311}
\ee
Conversely, for any covector $\xi \in T_{x}^{*}M\setminus \{0\}$, there exists a unique vector $y\in T_{x}M\setminus \{0\}$ such that $\xi =g_{y}(y, \cdot)\in T^{*}_{x}M$.
\end{lem}

Naturally, by Lemma \ref{shen311}, we define a map ${\cal L}: TM \rightarrow T^{*}M$ by
\[
{\cal L}(y):=\left\{
\begin{array}{ll}
g_{y}(y, \cdot), & y\neq 0, \\
0, & y=0.
\end{array} \right.
\]
It follows from (\ref{shenF311}) that
\be
F(x,y)=F^{*}(x, {\cal L}(y)).
\ee
Thus ${\cal L}$ is a norm-preserving transformation. We call ${\cal L}$ the Legendre transformation on Finsler manifold $(M, F)$.

Take a basis $\{{\bf b}_{i}\}^{n}_{i=1}$ for $TM$ and its dual basis $\{\theta ^{i}\}_{i=1}^{n}$ for $T^{*}M$. Express
\[
\xi ={\cal L}(y)=\xi _{i}\theta ^{i}=g_{ij}(x,y)y^{j}\theta ^{i},
\]
where $g_{ij}(x,y):=\frac{1}{2}\left[F^2\right]_{y^{i}y^{j}}(x,y)$. Let
\be
g^{*kl}(x,\xi):=\frac{1}{2}\left[F^{*2}\right]_{\xi _{k}\xi_{l}}(x,\xi).
\ee
For any $\xi ={\cal L}(y)$, differentiating $F^{2}(x,y)=F^{*2}(x,{\cal L}(y))$ with respect to $y^i$ yields
\[
\frac{1}{2}\left[F^2\right]_{y^{i}}(x,y)=\frac{1}{2}\left[F^{*2}\right]_{\xi _{k}}(x,\xi)g_{ik}(x,y), \label{shen3121}
\]
which implies
\be
g^{*kl}(x,\xi)\xi _{l}=\frac{1}{2}\left[F^{*2}\right]_{\xi _{k}}(x,\xi)=\frac{1}{2}g^{ik}(x,y)\left[F^2\right]_{y^{i}}(x,y)=y^{k}. \label{dualY}
\ee
Then,  we can get (see \cite{Cheng1}\cite{shen1})
\be
g^{*kl}(x,\xi)= g^{kl}(x,y), \label{Fdual}
\ee
where $(g^{kl}(x,y))= (g_{kl}(x,y))^{-1}.$

Let $(M, F, {m})$ be a Finsler manifold equipped with a measure ${m}$ on $M.$  Given an open set $\Omega \subset M,$ let $H_{\rm loc }^{1}(\Omega)$ be the space of weakly differentiable functions $u$ on $\Omega$ such that both $u$ and $F^{*}(d u)$ belong to $L_{\rm loc }^{2}(\Omega)$. We remark that $H_{\text {loc }}^{1}(\Omega)$ is  a linear space determined independent of the choices of $F$ and ${m}$ . Define the energy functional $\mathcal{E}_{\Omega}: H_{\rm loc }^{1}(\Omega) \longrightarrow[0, \infty]$ by
\be
\mathcal{E}_{\Omega}(u):=\frac{1}{2} \int_{\Omega} F^{*}(x, d u)^{2} d {m}.
\ee
Further, define the Sobolev space associated with $\mathcal{E}_{\Omega}$ by
\[
H^{1}(\Omega):=\left\{u \in L^{2}(\Omega) \cap H_{{\rm loc}}^{1}(\Omega) \mid {\cal E}_{\Omega}(u)+{\cal E}_{\Omega}(-u)<\infty\right\},
\]
and let $H_{0}^{1}(\Omega)$ be the closure of $\mathcal{C}_{c}^{\infty}(\Omega)$ with respect to the (absolutely homogeneous) Sobolev norm
$$
\|u\|_{H^{1}(\Omega)}:=\|u\|_{L^{2}(\Omega)}+\left\{\mathcal{E}_{\Omega}(u)+\mathcal{E}_{\Omega}(-u)\right\}^{1 / 2},
$$
where $\mathcal{C}_{c}^{\infty}(\Omega)$ denotes the set of $\mathcal{C}^{\infty}$-functions on $\Omega$ with compact support.

Given a smooth function $u$ on $M$, the differential $du_{x}$ at any point $x\in M$,
\[
du_{x}=\frac{\pa u}{\pa x^i}(x)dx^{i}
\]
is a linear function on $T_{x}M$. We define the gradient vector $\nabla u(x)$ of $u$ at $x\in M$ by $\nabla u(x):={\cal L}^{-1}\left(du(x)\right)\in T_{x}M$. In a local coordinate system, by (\ref{dualY}), we can express $\nabla u$ as
\be
\nabla u(x)=\left\{
\begin{array}{ll}
g^{*ij}(x,d u)\frac{\pa u}{\pa x^i}\frac{\pa}{\pa x^j}, & x\in M_{u},  \label{gradientV}\\
0, & x \in M\setminus M_{u},
\end{array}
\right.
\ee
where $M_{u}=\{x\in M \mid du (x)\neq 0\}.$ Further, by Lemma \ref{shen311}, we have the following
\be
du_{x}(v)=g_{\nabla u_{x}}(\nabla u_{x}, v), \ \ \ \forall v\in T_{x}M
\ee
and
\be
F(x, \nabla u_{x})=F^{*}(x, du_{x})=\frac{du_{x} (\nabla u_{x})}{F(x,\nabla u_{x})}.
\ee
We must be careful when $d u(x)=0,$ because $g^{*i j}(d u(x))$ is not defined and the Legendre transform $\mathcal{L}^{-1}$ is only continuous at the zero section. Besides, if $u \in \mathcal{C}^{l}(M),$ then ${\nabla}u$ is $\mathcal{C}^{l-1}$ on $M_{u}$ whereas only continuous on $M \backslash M_{u}$.

Associated with the measure $m$ on $M$, we first decomposed the measure $m$ as $d m = {\rm e}^{\Phi} dx^{1} dx^{2} \cdots d x^{n}$. Then the divergence of a differentiable vector field $V$ on $M$ is defined by
\be
{\rm div}_{m} V:=\sum_{i=1}^{n}\left(\frac{\partial V^{i}}{\partial x^{i}}+V^{i} \frac{\partial \Phi}{\partial x^{i}}\right), \quad V=\sum_{i=1}^{n} V^{i} \frac{\partial}{\partial x^{i}}. \label{divdef1}
\ee
One can also define ${\rm div}_{m} V$ in the weak form by following divergence formula:
\be
\int_{M} \phi \ {\rm div}_{m}V \ d{m}=-\int_{M} d\phi (V) \ d{m} \label{divdef2}
\ee
for all $\phi \in \mathcal{C}_{c}^{\infty}(M)$.

Now we define the Finsler Laplacian $\Delta u$ of $u \in H_{\rm loc}^{1}(M)$ by
\be
\Delta u:= {\rm div}_{m}(\nabla u). \label{laplace1}
\ee
Equivalently, we can define Laplacian $\Delta u$ on the whole $M$ in the weak sense by
\be
\int_{M} \phi \ {\Delta}u \ d{m}:=-\int_{M} d \phi({\nabla} u) d{m} \label{laplace2}
\ee
for all $\phi \in \mathcal{C}_{c}^{\infty}(M)$.

From (\ref{laplace1}), Finsler Laplacian is a nonlinear elliptic differential operator of the second order. Moreover, since the gradient vector field $\nabla u$ is merely continuous on $M \backslash M_{u}$, even when $u \in \mathcal{C}^{\infty}(M),$ it is necessary to introduce the Laplacian in the weak form as (\ref{laplace2}).

By the definitions, we have the following for any smooth function $\varphi$ on $M$,
\be
{\rm div}_{m}(\varphi \nabla u)=\varphi \Delta u+d\varphi (\nabla u). \label{DLOper}
\ee

Let $V=V^{i} \frac{\partial}{\partial x^{k}}$ be a nonzero measurable vector field on $M.$ One can introduce a weighted Riemannian metric $g_{V}$ on $M$ via
\[
g_{V}(X, Y)=g_{i j}(x, V) X^{i} Y^{j}, \quad \text{ for } X, Y \in T_{x} M.
\]
In particular, $g_{V}(V, V)=F^{2}(x, V)$.  For $u \in H_{\rm loc }^{1}(M)$ such that $du=0$ almost everywhere on $\{x \in M \mid V(x)=0\}$ (in other words, $V \neq 0$ almost everywhere on $M_{u}$ ), we can define the linearized gradient vector and the linearized Laplacian on the weighted Riemannian manifold $\left(M, g_{V}, {m}\right)$ by
\be
\nabla^{V} u:=\left\{\begin{array}{ll}\sum\limits_{i, j=1}^{n} g^{i j}(x, V) \frac{\partial u}{\partial x^{j}} \frac{\partial}{\partial x^{i}} & \text{on } M_{u}, \\
0 & \text{on } M \backslash M_{u},\end{array} \quad \Delta^{V} u:= {\rm div}_{m}\left(\nabla^{V} u\right)\right.. \label{weigraLa}
\ee

It is easy to get the following by (\ref{weigraLa}): for any $u \in H_{\rm loc}^{1}(M)$, we have
\be
\nabla^{\nabla u} u=\nabla u , \ \ \ \Delta^{\nabla u} u=\Delta u. \label{weiGra}
\ee
Further, for any $f_{1}, f_{2} \in H_{\mathrm{loc}}^{1}(M)$ satisfying $d f_{1}=d f_{2}=0$ almost everywhere on $M \backslash M_{u}$, we have
\be
d f_{2}\left(\nabla^{\nabla u} f_{1}\right)= d f_{1}\left(\nabla^{\nabla u} f_{2}\right). \label{LinGra}
\ee

\vskip 3mm

Finsler geometry is just Riemannian geometry without the quadratic restriction (\cite{Chern}). The Ricci curvature in Finsler geometry is just a natural extension of the Ricci curvature in Riemann geometry. However, we have a difficulty on the choice of a measure in Finsler geometry, because it is impossible to choose a unique canonical measure like the volume measure in the Riemannian setting. Naturally, by choosing an arbitrary measure $m$ on a Finsler manifold $(M,F)$, Ohta modified the Ricci curvature and defined the weighted Ricci curvature in Finsler geometry (\cite{Oh1}). Concretely, for an $n$-dimensional Finsler manifold $(M, F, m)$ equipped with a smooth measure $m$ and for any $v \in T_{x} M \backslash\{0\},$ let $\eta:(-\varepsilon, \varepsilon) \longrightarrow M$ be the geodesic with $\dot{\eta}(0)=v$ and decompose the measure ${m}$ along $\eta$ as
\[
m = {\rm e}^{-\psi_{\eta}} \sqrt{{\rm det}\left(g_{i j}(\eta , \dot{\eta})\right)}\ d x^{1} d x^{2} \cdots d x^{n},
\]
where $\psi_{\eta}= \psi_{\eta}\left(\eta (t), \dot{\eta}(t)\right): (-\varepsilon, \varepsilon) \longrightarrow R$ is a $\mathcal{C}^{\infty}$-function. Then, for $N \in R \backslash\{n\},$ define the weighted Ricci curvature
\be
{\rm Ric}_{N}(v):= {\rm Ric}(v)+\psi_{\eta}^{\prime \prime}(0)-\frac{\psi_{\eta}^{\prime}(0)^{2}}{N-n}. \label{weiRicci1}
\ee
As the limits of $N \rightarrow \infty$ and $N \downarrow n,$  we define the weighted Ricci curvatures as follows.
\be
\begin{aligned}
{\rm Ric}_{\infty}(v) &: = {\rm Ric}(v)+\psi_{\eta}^{\prime \prime}(0), \\
{\rm Ric}_{n}(v)& :=\left\{\begin{array}{ll} {\rm Ric}(v)+\psi_{\eta}^{\prime \prime}(0) & \text { if } \psi_{\eta}^{\prime}(0)=0, \\
-\infty & \text { if } \psi_{\eta}^{\prime}(0) \neq 0. \end{array}\right.
\end{aligned}
\ee

It should be point out that the quantity ${\bf S}(x, v):=\psi_{\eta}^{\prime}(0)$ is just the S-curvature with respect to the measure $m$ in Finsler geometry(\cite{ChernShen}). We say that ${\rm Ric}_{N} \geq K$ for $K \in {R}$ if ${\rm Ric}_{N}(v) \geq K F^{2}(x, v)$ for all $x\in M$ and $v \in T_{x} M$.

\section{Poincar\'{e}-Lichnerowicz inequality}\label{Sec3}

Let $(M, F, m)$ be a Finsler manifold equipped with a measure $m$ on $M$. In order to prove Theorem \ref{PLineq}, we need some necessary lemmas. Firstly, we have the following lemma.

\begin{lem}{\rm (\cite{Oh4}\cite{OS1})} \label{Lemma11.4} If $(M, F, m)$  satisfies that $\Lambda_{F}<\infty, \ {m}(M)< \infty$ and is complete, then the constant function $1$ belongs to $H_{0}^{1}(M)$, where $\Lambda_{F}$ denotes the reversibility constant of $F$.
\end{lem}

Further,  we have the following result.

\begin{lem}\label{fsquareL}
For any $f, u \in H_{\mathrm{loc}}^{1}(M)$ satisfying $d f= 0$ almost everywhere on $M \backslash M_{u}$, we have
\be
\Delta^{\nabla u}f^{2}=2f \Delta^{\nabla u}f + 2 g_{\nabla u}\left(\nabla ^{\nabla u}f ~, \nabla ^{\nabla u}f\right), \label{eqLemma3.1}
\ee
equivalently,
\be
f \Delta^{\nabla u}f = \frac{1}{2}\Delta^{\nabla u}f^{2}- g_{\nabla u}\left(\nabla ^{\nabla u}f ~, \nabla ^{\nabla u}f\right). \label{squareLa}
\ee
\end{lem}

\noindent{\it Proof.} By (\ref{weigraLa}), we have
\[
\nabla^{\nabla u}f^{2} = 2f \nabla^{\nabla u}f.
\]
Further, by (\ref{DLOper}) and (\ref{weigraLa}), we have
\beqn
\Delta ^{\nabla u} f^{2}&=& {\rm div}_{m}\left(\nabla^{\nabla u}f^{2}\right)= {\rm div}_{m}(2 f \nabla^{\nabla u}f)\\
&=& 2f \Delta ^{\nabla u} f + 2 df (\nabla^{\nabla u}f).
\eeqn
It is easy to see that
\[
 df (\nabla^{\nabla u}f)=g_{\nabla u}\left(\nabla ^{\nabla u}f ~, \nabla ^{\nabla u}f\right).
\]
Thus we get (\ref{eqLemma3.1}). \qed

\vskip 2mm

The following lemma is extremely important for our discussions which is derived from Lemma \ref{IntformBo}.
\begin{lem}
Assume that $(M, F, {m})$ is compact and satisfies ${\rm Ric}_{\infty} \geq K>0$. Then , for any $u \in H^{2}(M) \cap \mathcal{C}^{1}(M)$ such that $\Delta u \in H^{1}(M)$, we have
\be
\int_{M} F^{2}({\nabla}u) d{m} \leq \frac{1}{K}\left\{ \int_{M}({\Delta} u)^{2} d{m} - \int_{M} g_{\nabla u}\left(\nabla ^{\nabla u}F(\nabla u) , \nabla ^{\nabla u}F(\nabla u)\right)dm \right\}. \label{mainlem}
\ee
\end{lem}

\noindent{\it Proof.} From (\ref{Intform}) and by taking test function $\phi \equiv 1$ according to Lemma \ref{Lemma11.4}, we have the following
\be
\int_{M} \left\{d({\Delta} u)({\nabla} u)+K F^{2}({\nabla}u)+d[F({\nabla}u)]\left(\nabla^{\nabla u}[F({\nabla} u)]\right)\right\} d{m} \leq 0.\label{firstequ}
\ee
Noticing that
\beqn
&& \int_{M} d({\Delta} u)({\nabla} u)dm = - \int_{M} (\Delta u)^{2}dm = -\|\Delta u \|^{2}_{L^{2}},\\
&& \int_{M} d[F({\nabla}u)]\left(\nabla^{\nabla u}[F({\nabla} u)]\right) = - \int_{M} F({\nabla}u) \Delta^{\nabla u}\left[F({\nabla}u)\right]dm ,
\eeqn
(\ref{firstequ}) becomes
\[
 K \int_{M}F^{2}({\nabla}u)dm \leq  \int_{M} (\Delta u)^{2}dm +   \int_{M} F({\nabla}u) \Delta^{\nabla u}\left[F({\nabla}u)\right]dm .
\]

Further, by Lemma \ref{fsquareL}, we have
\beqn
K \int_{M}F^{2}({\nabla}u)dm &\leq&  \int_{M} (\Delta u)^{2}dm +   \int_{M} \Delta ^{\nabla u}\left[\frac{ F^{2}({\nabla}u)}{2} \right]dm \\
    && -\int_{M} g_{\nabla u}\left(\nabla ^{\nabla u}F(\nabla u) , \nabla ^{\nabla u}F(\nabla u)\right)dm \\
   &=& \int_{M} (\Delta u)^{2}dm  -\int_{M} g_{\nabla u}\left(\nabla ^{\nabla u}F(\nabla u) , \nabla ^{\nabla u}F(\nabla u)\right)dm .
\eeqn
From this, we obtain (\ref{mainlem}). \qed

\vskip 2mm

In the following, normalizing $m$ as $m(M)=1,$ we define the variance of $f \in L^{2}(M)$ as
\[
\operatorname{Var}_{m}(f):=\int_{M}\left(f-\int_{M} f d m \right)^{2} d m = \int_{M} f^{2} d m -\left(\int_{M} f d m \right)^{2}.
\]
At the same time, we will mainly consider the outcome of a kind of ergodicity
\be
\lim\limits_{t\rightarrow \infty} u_{t} = \int_{M} u_{0} d{m} \ \ \ {\rm in} \ L^{2}(M) \label{eq15.2}
\ee
for all global solutions $\left(u_{t}\right)_{t \geq 0}$ to the heat equation $\pa _{t}u_{t} = \Delta u_{t}$ (see \cite{Oh2}).

Now we are in the position to prove Theorem \ref{PLineq}.
\vskip 2mm

\noindent {\it Proof of Theorem \ref{PLineq}.} \ Let $(u_{t})_{t\geq 0}$ be the global solution to the heat equation with $u_{0}= f$. Put $\Phi (t):= \|u_{t}\|^{2}_{L^2}= \int_{M} u^{2}_{t} d m$. Then the ergodicity (\ref{eq15.2}) implies that
\beqn
{\rm Var}_{m}(f)&=&  \int_{M} f^{2} d m -\left(\int_{M} f d m \right)^{2}\\
                &=&  \Phi (0) - \lim\limits_{t \rightarrow \infty} \Phi (t)= - \int_{0}^{\infty} \Phi ' (t) dt .
\eeqn
On the other hand, by the definition of $\Phi (t)$ and (\ref{squareLa}), we know that
\be
\Phi^{\prime}(t)=2 \int_{M} u_{t} \Delta u_{t} d{m}=-2 \int_{M} F^{2}\left(\boldsymbol{\nabla} u_{t}\right) d{m}=-4 \mathcal{E}\left(u_{t}\right).\label{Phiprime}
\ee

By (4.2) in \cite{OS2}, we have the following
\be
\frac{\partial}{\partial t}\left[F^{2}\left({\nabla} u_{t}\right)\right]=2 d\left({\Delta} u_{t}\right)\left({\nabla} u_{t}\right) \label{Le14.1}
\ee
for all $t >0$. Hence, from (\ref{Phiprime}) and by (\ref{Le14.1}), we have
\[
\Phi^{\prime \prime}(t)=-4 \int_{M} d\left(\Delta u_{t}\right)\left(\nabla u_{t}\right) d {m}=4 \int_{M}\left(\Delta u_{t}\right)^{2} dm =4\left\|\Delta u_{t}\right\|_{L^{2}}^{2}
\]
for all $t>0$. Thus, from (\ref{mainlem}) for $u_{t}$ and by (\ref{Phiprime}), we get the following
\[
-\frac{1}{2}\Phi^{\prime}(t) \leq \frac{1}{4K} \Phi^{\prime \prime}(t)-\frac{1}{K}\int_{M}g_{\nabla u_{t}}\left(\nabla ^{{\nabla u_{t}}}F({\nabla u_{t}}), \nabla ^{{\nabla u_{t}}}F({\nabla u_{t}})\right)dm,
\]
that is,
\[
- \Phi^{\prime}(t) \leq \frac{1}{2K} \Phi^{\prime \prime}(t)-\frac{2}{K}\int_{M}g_{\nabla u_{t}}\left(\nabla ^{{\nabla u_{t}}}F({\nabla u_{t}}), \nabla ^{{\nabla u_{t}}}F({\nabla u_{t}})\right)dm.
\]
Then,
\beqn
&& {\rm Var}_{m}(f)=-\int_{0}^{\infty} \Phi ^{\prime}(t)dt \\
&& \leq \frac{1}{2K}\left(\lim\limits_{t \rightarrow \infty}\Phi ^{\prime}(t) - \Phi ^{\prime}(0)\right)-\frac{2}{K}\int_{M}\left(\int_{0}^{\infty} g(t)dt\right)dm,
\eeqn
where
\[
g(t):= g_{\nabla u_{t}}\left(\nabla ^{{\nabla u_{t}}}F({\nabla u_{t}}), \nabla ^{{\nabla u_{t}}}F({\nabla u_{t}})\right).
\]
From (\ref{Phiprime}), $\Phi ^{\prime}(0)= -4 {\cal E}(f)$. Further, it is not difficult to prove that $\lim\limits_{t\rightarrow \infty}{\cal E}(u_{t})=0$ (e.g. see Proposition 13.15 in \cite{Oh4}).  Thus, we get the following
\[
{\rm Var}_{m}(f)\leq \frac{2}{K}{\cal E}(f)- \frac{2}{K}\int_{M}\left(\int_{0}^{\infty} g(t)dt\right)dm.
\]
This completes the proof of Theorem \ref{PLineq}. \qed
\vskip 2mm
In \cite{Oh2}, Ohta obtained Poincar\'{e}-Lichnerowicz inequality under the condition that ${\rm Ric}_{N} \geq K>0$ for some $N \in(-\infty, 0) \cup[n, \infty]$ as follows
\be
{\rm Var}_{m}(f) \leq \frac{N-1}{K N} \int_{M} F^{2}({\nabla} f) d{m}. \label{eq15.6}
\ee
In the case of $N=\infty,$ the inequality (\ref{eq15.6}) becomes as
\be
{\rm Var}_{m}(f) \leq \frac{1}{K} \int_{M} F^{2}({\nabla} f) d{m}. \label{eq15.8}
\ee
Notice that (\ref{PoLiin}) is stronger than (\ref{eq15.8}).

\section{Logarithmic Sobolev inequality}\label{Sec4}

In this section we will study the logarithmic Sobolev inequality under the condition that ${\rm Ric}_{\infty} \geq K >0 $. The logarithmic Sobolev inequality plays a quite important role in infinite-dimensional analysis and probability theory. From the geometric viewpoint, the logarithmic Sobolev inequality plays an important role in the study of entropy in Ricci flow theory. For brevity, we introduce the following notation (see \cite{Oh2}):
\be
\Gamma _{2}(u):= \Delta^{\nabla u}\left[\frac{F^{2}({\nabla} u)}{2}\right]-d({\Delta} u)({\nabla} u). \label{Gamma2}
\ee
Then the improved Bochner inequality (\ref{eq12.13}) can be rewritten as follows:
\be
\Gamma _{2}(u) \geq K F^{2}({\nabla} u)+d[F({\nabla} u)]\left(\nabla^{\nabla u}[F({\nabla} u)]\right). \label{eq12.132}
\ee
Further, for nonnegative function $f \in L^{1}(M)$ with $\int_{M} f d{m}=1,$  define the relative entropy with respect to ${m}$ as
\be
{\rm Ent}_{m}(f{m}):=\int_{M} f \log f d{m}. \label{reentropy}
\ee

\begin{lem}\label{Prop16.3} {\rm (\cite{Oh2}\cite{Oh4})} \ Assume that $(M, F, {m})$ is compact and satisfies ${m}(M)=1$ and
\be
\int_{M} \frac{F^{2}({\nabla} u)}{u} d{m} \leq -C \int_{M}\left\{d u\left(\nabla^{\nabla u}\left[\frac{F^{2}({\nabla} u)}{2 u^{2}}\right]\right)+u \cdot d[{\Delta}(\log u)]({\nabla}[\log u])\right\} d {m} \label{eq16.3}
\ee
for some constant $C>0$ and all positive functions $u \in H^{2}(M) \cap \mathcal{C}^{1}(M)$ such that $\Delta u \in H^{1}(M)$. Then the logarithmic Sobolev inequality
\be
\int_{M} f \log f d {m} \leq \frac{C}{2} \int_{M} \frac{F^{2}({\nabla} f)}{f} d{m} \label{eq16.4}
\ee
holds for all nonnegative functions $f \in H^{1}(M)$ with $\int_{M} f d {m}=1$.
\end{lem}

By the notation (\ref{Gamma2}), the condition (\ref{eq16.3}) is equivalent to the following
\be
\int_{M} u F^{2}(\nabla[\log u]) d{m} \leq C \int_{M} u \Gamma_{2}(\log u) d{m}. \label{eq16.5}
\ee

In \cite{Oh2}, Ohta proved the following logarithmic Sobolev inequality under the condition that ${\rm Ric}_{N} \geq K>0$ for some $N \in[n, \infty)$ and ${m}(M)=1$:
$$
\int_{\{f>0\}} f \log f d{m} \leq \frac{N-1}{2 K N} \int_{\{f>0\}} \frac{F^{2}({\nabla} f)}{f} d{m}
$$
for all nonnegative functions $f \in H^{1}(M)$ with $\int_{M} f d{m}=1$. In the case of $N=\infty,$ the logarithmic Sobolev inequality can be rewritten as folows
\be
\int_{M} f \log f d{m} \leq \frac{1}{2 K} \int_{M} \frac{F^{2}({\nabla} f)}{f} d{m}. \label{th18.8}
\ee
In fact, Ohta has given a direct proof of (\ref{th18.8}) by the curvature-dimension condition in \cite{Oh1}.

In the following, we will give a new and different proof of (\ref{th18.8}) (that is, (\ref{eq16.7m})) by using the improved Bochner inequality (\ref{eq12.132}).
\vskip 3mm
\noindent{\it Proof of Theorem \ref{LSIne}.} \ Fix $h \in \mathcal{C}^{\infty}(M)$ and consider the function ${\rm e}^{a h}$ for $a>0$. It is easy to see that
\be
{\nabla}\left({\rm e}^{a h}\right)=a {\rm e}^{a h} {\nabla} h, \quad  {\Delta}\left({\rm e}^{a h}\right)=a {\rm e}^{a h}\left({\Delta} h+aF^{2}({\nabla} h)\right). \label{eq16.8}
\ee
Then, by (\ref{Gamma2}), we have
\beq
\Gamma_{2}\left({\rm e}^{a h}\right)&=&\Delta^{\nabla h}\left[\frac{a^{2} {\rm e}^{2 a h} F^{2}({\nabla} h)}{2}\right]-a^{2} {\rm e}^{a h} d\left[{\rm e}^{a h}\left(\Delta h+a F^{2}({\nabla} h)\right)\right]({\nabla} h) \nonumber\\
&=& a^{2} {\rm e}^{2 a h}\left\{\Gamma_{2}(h)+a d\left[F^{2}({\nabla} h)\right]({\nabla} h)+a^{2} F^{4}({\nabla} h)\right\}. \label{eq16.9}
\eeq
Now, applying the Bochner inequality (\ref{eq12.132}) to $u= {\rm e}^{ah}$ yields
\beqn
&& a^{2} {\rm e}^{2 a h}\left\{\Gamma_{2}(h)+a d\left[F^{2}({\nabla} h)\right]({\nabla} h)+a^{2} F^{4}({\nabla} h)\right\} \\
&& \geq K a^{2} {\rm e}^{2 a h}F^{2}({\nabla} h)+ a d\left[{\rm e}^{a h}F({\nabla} h)\right]\left(\nabla ^{\nabla h}\left[a {\rm e}^{a h}F({\nabla} h)\right]\right) \\
&&= K a^{2} {\rm e}^{2 a h}F^{2}({\nabla} h)+ a^{2} {\rm e}^{2 a h}\left[a F({\nabla} h) dh + dF(\nabla h)\right]\left(a  F({\nabla} h) \nabla h + \nabla^{\nabla h}F(\nabla h)\right) \\
&&= a^{2} {\rm e}^{2 a h} \left\{K F^{2}({\nabla} h)+ a^{2}F^{4}({\nabla} h)+ 2 a F({\nabla} h)d[F({\nabla} h)](\nabla h) \right. \\
&& \ \ \ \ \left. + d[F({\nabla} h)]\left(\nabla^{\nabla h}F(\nabla h)\right)  \right\} \\
&&= a^{2} {\rm e}^{2 a h} \left\{K F^{2}({\nabla} h)+ a^{2}F^{4}({\nabla} h)+ a d[F^{2}({\nabla} h)](\nabla h) \right. \\
&& \ \ \ \ \left.+ g_{\nabla h}\left(\nabla^{\nabla h}F(\nabla h)~, \nabla^{\nabla h}F(\nabla h) \right)\right\},
\eeqn
where we have used the following facts
\[
dh\left(\nabla^{\nabla h}F(\nabla h)\right)= d[F({\nabla} h)](\nabla h),
\]
\[
d[F({\nabla} h)]\left(\nabla^{\nabla h}F(\nabla h)\right)= g_{\nabla h}\left(\nabla^{\nabla h}F(\nabla h)~, \nabla^{\nabla h}F(\nabla h) \right).
\]
Then we obtain the inequality
\be
\Gamma _{2}(h)\geq K F^{2}({\nabla} h)+ g_{\nabla h}\left(\nabla^{\nabla h}F(\nabla h)~, \nabla^{\nabla h}F(\nabla h) \right). \label{eq16.11}
\ee
It is surprising that none of the terms with $a$ occur in (\ref{eq16.11}). Further, we obtain
\beqn
\int_{M} {\rm e}^{h} \Gamma_{2}(h) d{m} &\geq & K \int_{M} {\rm e}^{h} F^{2}({\nabla} h) d{m}+ \int_{M} {\rm e}^{h}g_{\nabla h}\left(\nabla^{\nabla h}F(\nabla h)~, \nabla^{\nabla h}F(\nabla h) \right)dm \\
&\geq & K \int_{M} {\rm e}^{h} F^{2}({\nabla} h) d{m}.
\eeqn
This is just the inequality (\ref{eq16.5}) with $C= 1 / K$ for $u= {\rm e}^{h}$.  By approximation this implies (\ref{eq16.5}) for all positive function $u$ with the required properties. This completes the proof by Lemma \ref{Prop16.3}. \qed

\section{The estimate of the volumes of geodesic spheres}\label{Sec5}

Let $(M, F, m)$ be a Finsler manifold equipped with a measure $m$ on $M$.  Fix $p \in M$. For a unit vector $v \in I_{p} M := \{y\in T_{p}M \mid F(p, y)=1\}$, define
$$
c(v):=\sup \left\{t>0 \mid d_{F}\left(p, \exp _{p}(t v)\right)=t\right\} \in(0, \infty]
$$
If $c(v)<\infty,$ we call the point $\exp _{p}(c(v) v)$ a cut point of $p.$  The set of all cut points of $p$ is called the cut locus of $p$ and will be denoted by ${\rm Cut}(p),$ that is
\[
{\rm Cut}(p):=\{\exp _{p}(c(v)v) \mid v\in I_{p} \}.
\]
The injectivity radius $i_{p}$ at $p$ is defined by $i_{p}=\inf \left\{c(v) \mid v \in I_{p}\right\}$.

For $R> 0$, put
\[
B_{p}(R):= \exp _{p}\left[y\in T_{p}M \mid F(p, y)\leq R\right].
\]
$B_{p}(R)$ is called the geodesic ball of radius $R$ at the center $p$. Obviously, geodesic ball $B_{p}(R)$ is well-defined when $R < i_{p}$. In this paper, the radius $R$  of  any geodesic ball $B_{p}(R)$ always satisfies that $R < i_{p}$.

Let $r(x):= d_{F}(p, x)$ be the distance function on $M$ from $p$. By Lemma 3.2.3 in \cite{shen1}, we know that $F(x, \nabla r(x))=1$. Starting from this point, we can prove Theorem \ref{volgeoBall}.

\vskip 3mm

\noindent{\it Proof of Theorem \ref{volgeoBall}.} By Lemma \ref{IntformBo},  taking test function $\phi \equiv 1$ in (\ref{Intform}) yields the following
\[
\int_{M} \left\{d({\Delta} u)({\nabla} u)+K F^{2}({\nabla}u)+d[F({\nabla}u)]\left(\nabla^{\nabla u}[F({\nabla} u)]\right)\right\} d{m} \leq 0.
\]
Now we employ the distance function $r(x)= d_{F}(p, x)$ as $u$ in the above inequality and restrict our discussions on the geodesic ball $B_{p}(R)$. Then,  by the fact that $F(x, \nabla r(x))=1$, we have the following
\[
\int_{B_{p}(R)}\left\{d({\Delta} r)({\nabla} r)+K \right\}dm \leq 0,
\]
That is,
\[
-\int_{B_{p}(R)}\left(\Delta r\right)^{2} dm+ K {\rm vol}\left(B_{p}(R)\right)\leq 0.
\]
Then we get
\[
{\rm vol}\left(B_{p}(R)\right)\leq \frac{1}{K}\int_{B_{p}(R)}\left(\Delta r\right)^{2} dm .
\]
This is just (\ref{geoBall}). \qed

\vskip 8mm

\noindent {\bf Acknowledgements.} The author would like to thank Professor Shin-ichi Ohta for his valuable suggestions.

\vskip 16mm

\noindent
Xinyue Cheng \\
School of Mathematical Sciences \\
Chongqing Normal University \\
Chongqing  401331,  P. R. of China  \\
E-mail: chengxy@cqnu.edu.cn


\vskip 16mm
\begin{thebibliography}{Ma}


\bibitem{Cheng1} X, Cheng, Some fundamental problems in global Finsler geometry, arXiv: 1910.08267v1 [math.DG], 18 Oct 2019.


\bibitem{Chern} S. S. Chern, Finsler geometry is just Riemannian geometry without the quadratic restriction, Notices of the American Mathematical Society, September 1996.

\bibitem{ChernShen} S. S. Chern and Z. Shen, Riemann-Finsler Geometry, Nankai Tracts in Mathematics, Vol. 6, World Scientific, 2005.

\bibitem{Oh1} S. Ohta, Finsler interpolation inequalities, Calc. Var. Partial Differential Equations, {\bf 36}(2009), 211-249.

\bibitem{Oh2} S. Ohta, Some functional inequalities on non-reversible Finsler manifolds, Proc. Indian Acad. Sci. Math. Sci., {\bf 127}(2017), 833-855.

\bibitem{Oh3} S. Ohta,  A semigroup approach to Finsler geometry: Bakry-Ledoux's isoperimetric inequality, Comm. Anal. Geom., to appear.

\bibitem{Oh4} S. Ohta, Comparison Finsler geometry, in preparation.

\bibitem{OS1} S. Ohta and K.-T. Sturm, Heat flow on Finsler manifolds, Comm. Pure Appl. Math., {\bf 62}(2009), 1386-1433.

\bibitem{OS2} S. Ohta and K.-T. Sturm, Bochner-Weitzenb$\ddot{o}$ck formula and Li-Yau estimates on Finsler manifolds, Advances in Mathematics, {\bf 252}(2014), 429-448.

\bibitem{shen1} Z. Shen, Lectures on Finsler Geometry, World Scientific, Singapore, 2001.

\bibitem{WangXia} G. Wang and C. Xia, A sharp lower bound for the first eigenvalue on Finsler manifolds, Ann. I. H. Poincar\'{e}-Nonlinear Anal., {\bf 30}(2013), 983 -996.

\bibitem{Xia} Q. Xia, Geometric and functional inequalities on Finsler manifolds, The Journal of Geometric Analysis, {\bf 30}(2020), 3099 - 3148.

\end{thebibliography}
\end{document}